\documentclass[12pt]{amsart}
\usepackage{amssymb}
\usepackage{amsthm}
\usepackage{euscript}
\usepackage{xspace}
\usepackage{enumerate}
\def\rmark{\mbox{$\rm\bf\rule{0.06em}{1.45ex}\kern-0.05em R$}}
\def\nmark{\mbox{$\rm\bf\rule{0.06em}{1.45ex}\kern-0.05em N$}}
\def\hmark{\mbox{$\rm\bf\rule{0.06em}{1.45ex}\kern-0.05em H$}}
\def\cmark{\mbox{$\rm\bf\kern0.2em\rule{0.06em}{1.45ex}\kern-0.3em C$}}

\theoremstyle{plain}
\newtheorem{lem}[subsection]{Lemma}

\newtheorem{thm}[subsection]{Theorem}
\newtheorem{prop}[subsection]{Proposition}

\newtheorem{cor}[subsection]{Corollary}

\begin{document}
\title[nest modules of matrices]{On nest modules of matrices over division rings}
\author{M. Rahimi-Alangi   and Bamdad R. Yahaghi}
\address{ Department of Mathematics, Payame Noor University, P.O. Box 19395-3697 Tehran, Iran \newline
Department of Mathematics, Faculty of Sciences, Golestan University, Gorgan 19395-5746, Iran}
\email{mrahimi40@yahoo.com, \newline bamdad5@hotmail.com, bamdad@bamdadyahaghi.com}

\dedicatory{
``Heydar Baba may the sun warm your back,  \\
Make your smiles and your springs shed tears,\\
Your children collect a bunch of flowers, \\
Send it with the coming wind towards us, \\
Perhaps my sleeping furtune would awaken!\\
...\\
Heydar Baba may you be fortunate! \\
 Be surrounded with springs and orchards!\\
 May you live long after us!"\\
-Shahriar\\
With kind regards,\\ dedicated to Heydar Radjavi on the occasion of his eightieth birthday}

\keywords{Bimodule of rectangular matrices over a division ring, (Left/Right) Submodule,  Subbimodule, (One-sided) Ideal, Nest modules.\\
``Heydar Baba" is the name of a mountain overlooking the village Khoshgenab near Tabriz, where the well-known Iranian poet M.H. Shahriar was born and grew up. ``Heydar Baba Salam" is  the title of one of the most famous poems by Shahriar in Azeri Turkish in which he remembers his childhood and his memories from the mountain Heydar Baba and the village Khoshgenab. The translation is taken from the website of the Department of Near Eastern Studies of University of Michigan, Ann Arbor.  In Persian ``Baba'' means father. ``Baba'' is also an honorific term to address and refer to Sufi saints and mystics, e.g., Baba Taher. 
 }

\bibliographystyle{plain}

\begin{abstract}
Let $ m , n \in \mathbb{N}$, $D$ be a division ring,  and $M_{m \times n}(D)$ denote the bimodule of all  $m \times n$ matrices with entries from $D$. First, we characterize one-sided submodules of  $M_{m \times n}(D)$ in terms of left row reduced echelon or right column reduced echelon matrices with entries from $D$. Next, we introduce the notion of a nest module of matrices with entries from $D$. We then characterize submodules of nest modules of matrices over $D$ in terms of certain finite sequences of  left row reduced echelon or right column reduced echelon matrices with entries from $D$. We use this result to characterize principal submodules of nest modules. We also describe subbimodules of nest modules of matrices. As a consequence, we characterize (one-sided) ideals of nest algebras of matrices over division rings.
\end{abstract}
\maketitle \vspace{.5cm}

\bigskip

\begin{section}
{\bf One-sided submodules of $M_{m \times n}(D)$}
\end{section}

In this paper, we consider one-sided  submodules of $M_{m \times n}(D)$, where $D$ is a division ring. First, we present a characterization of one-sided submodules of $M_{m \times n}(D)$ via left row reduced or right column reduced echelon matrices with entries from $D$ (Theorem \ref{1.2}). Then we introduce nest modules of matrices and provide a characterization of their one-sided and two-sided submodules (Theorems \ref{2.2} and \ref{3.2}). As a consequence of our results, we characterize principal one-sided submodules of nest modules of matrices and in particular principal one-sided ideals of nest algebras of matrices (Theorem \ref{2.6}). It turns out that sub-bimodules of nest modules of matrices and in particular two-sided ideals of nest algebras of matrices are principal (Theorem \ref{3.2}). We have made a reasonably thorough search of the existing results, but to the best of our knowledge our results are new. For related results,  see \cite{Mi}, \cite{S}, and \cite{MSvW}. Through our results, one  sees that the set of all left row reduced (resp. right column reduced) echelon matrices over a division ring  forms a modular lattice via the operations join and meet which are defined in view of Theorem \ref{1.1} below. It seems that the lattice structure of the set of all left row reduced (resp. right column reduced) echelon matrices had not been noticed before.  

Let us set the stage by establishing some notation and definitions. Throughout this note, unless otherwise stated, $D$ denotes a
division ring, $F = Z(D)$ stands for the center of $D$, $m , n \in \mathbb{N}$,  $ 1 \leq i \leq m , 1 \leq j \leq n$, $M_{m \times n}(D)$, or simply $M_{m \times n}$,  denotes the set of all $m \times n$ matrices with entries from $D$, and $M_n(D) := M_{n \times n}(D)$. If there is no fear of confusion, we omit $D$ in our notations for the sake of simplicity. We view  $M_{m \times n}$ as an ($M_m,  M_n$)-bimodule via the matrix multiplication. In particular,  $M_n$ is an ($M_n,  M_n$)-bimodule, or simply an $M_n$-bimodule via the matrix multiplication. Also, in particular,  $D^n := M_{n \times 1}(D)$ and $ D_n:= M_{1 \times n} (D)$ are, respectively,  viewed as right and left $D$-modules, in other words right and left  vector spaces over $D$. As is usual, $E_{ij} \in M_{m \times n}$ denotes the matrix with $1$ in the $(i, j)$ place and zero elsewhere.  We use $ 0_{ m \times n}$ or $0_{mn} $ to denote the zero matrix in $M_{m \times n}$. We use $I_{m \times n}$ or $I_{mn}$ to denote the $m \times n$ matrix with $ 1$ in its $(i, i)$ place for each $ 1 \leq i \leq \min (m, n)$ and zero elsewhere. By convention, $ 0_n := 0_{ n \times n}$ and $ I_n := I_{n \times n}$, which is the identity matrix in $M_n$.   We call $E_{ij}$ a standard matrix. Also, by writing $E_{ij} \in M_n$  we clearly assume that $1 \leq i, j \leq n$. If $ A \in  M_{m \times n}$, we use
$ {\rm row}_i (A)$ and $ {\rm col}_j(A) $ to, respectively, denote the $i$th row and the $j$th column of $A$. For $ \mathcal{F} \subseteq M_{m \times n}$, by definition
$$ {\rm row}_i (\mathcal{F}) := \Big\{  {\rm row}_i (A):  A \in \mathcal{F}\Big\} , \  {\rm col}_j(\mathcal{F}) := \Big\{  {\rm col}_j(A): A \in \mathcal{F}\Big\}.$$

It is easily verified that for $ A  \in M_{m \times n}$ and $ E_{ij} \in M_{n}$ (resp. $ E_{ij} \in M_{m}$),  $A E_{ij}$ (resp. $ E_{ij}A$)
is the matrix whose $j$th column (resp. $i$th row) is the $i$th column (resp. $j$th row) of $A$ and every other column (resp. row) of it is zero, i.e., the operation $A  \rightarrow A E_{ij}$ (resp. $A  \rightarrow E_{ij}A$) takes the  $i$th column (resp. $j$th row) of $A$ to  the $j$th column (resp. $i$th row) of $A$ and takes every other column (resp. row) of $A$ to zero. Before stating our main result in this section, we make an easy  observation. Let $ \mathcal A $ and $\mathcal A'$ be subsets of $ M_{m \times n}$ that  are closed with respect to the addition of matrices and that absorb the multiplication by $E_{jj}$'s (resp. $E_{ii}$'s) from the right (resp. from the left), i.e.,  $ \mathcal {A} E_{jj} \subseteq \mathcal{A}$ and $\mathcal {A}' E_{jj} \subseteq \mathcal {A}'$ (resp. $ E_{ii} \mathcal{A} \subseteq \mathcal{A}$ and $  E_{ii} \mathcal{A}' \subseteq \mathcal{A}'$) for each $ 1 \leq j \leq n$ (resp. $ 1 \leq i \leq m$). Then $\mathcal{A} = \mathcal{A}'$ if and only if  ${\rm col}_j(\mathcal{A}) = {\rm col}_j(\mathcal{A}')$  for each $ 1 \leq j \leq n$, or $ {\rm row}_i (\mathcal{A}) = {\rm row}_i (\mathcal{A}')$ for each $ 1 \leq i \leq m$.

We use the symbols $ \mathcal{R}_{m \times n}(D)$ or simply $ \mathcal{R}_{m \times n}$ (resp. $ \mathcal{C}_{m \times n}(D)$ or simply $ \mathcal{C}_{m \times n}$) to denote the set of all $ m \times n$  left row reduced (resp. right column reduced) echelon matrices with entries from $D$.  Also, we write $ \mathcal{R}_{ n}:=\mathcal{R}_{n \times n}$ and
$ \mathcal{C}_{ n}:=\mathcal{C}_{n \times n}$. For $A \in M_{m \times n} $, we use $ {\rm LRS} (A)$ and $ {\rm RCS}(A)$  to, respectively,  denote the left row space and the right column space of $A$, i.e., the space spanned by the rows (resp. columns) of $A$. The matrices $ A, B \in M_{m \times n}$ are said to be left row (resp. right column) equivalent, and we write $A =_{lr} B $ (resp. $A=_{rc} B$), if $ {\rm LRS} (A)= {\rm LRS} (B)$ (resp.  $ {\rm RCS}(A)=  {\rm RCS}(B)$).   For $A, B \in  M_{m \times n}$, we write  $ A \leq_{lr} B $ if the left row space of $A$ is contained in the left row space of $B$. One can define $A \leq_{rc} B$ in a similar fashion. Clearly, $ A =_{lr} B $ (resp. $A=_{rc} B$) if and only if $ A \leq_{lr} B $  and $ B \leq_{lr} A $ (resp. $A \leq_{rc} B$ and $B \leq_{rc} A$).

Let $ A \in M_{m \times n} $ and $ B \in M_{p \times n}$.  We leave it as an exercise to the interested reader to show that
$ {\rm LRS} (A) \subseteq  {\rm LRS} (B)$ if and only if $ A = CB $ for some $ C \in M_{m \times p} $. Also if $  A \in M_{m \times n} $ and $ B \in M_{m \times p}$, then  $ {\rm RCS}(A) \subseteq   {\rm RCS}(B)$ if and only if $A = BC$ for some $ C \in  M_{p \times n}  $. For slight generalizations of infinite-dimensional counterparts of these facts in the setting of linear transformations see \cite[Corollary 1.4]{RY2}.  In particular, if $ A, B \in M_{m \times n}  $, then
$ A \leq_{lr} B $ (resp. $A \leq_{rc} B$) if and only if $ A = CB $ (resp. $A = BC$) for some $ C \in M_m$ (resp. $C \in M_n$).

We need the following well-known theorem for our main result in this section. We present a proof for reader's convenience. For the counterpart of the theorem below over general fields see  \cite[Theorem 2.5.11]{HK}. The first proof of the uniqueness in the following theorem is essentially taken from \cite{M}.

\bigskip

\begin{thm}\label{1.1}
{\rm  (i)}  {\it Let $W $ be a left subspace of $D_n$ with $ \dim W \leq m$. Then there exists a unique  left row reduced echelon matrix $ R \in M_{m \times n}(D)$ such that $ {\rm LRS}(R) = W$. }

{\rm (ii)}  {\it Let $W $ be a right subspace of $D^m$ with $ \dim W \leq n$. Then there exists a unique  right column reduced echelon matrix $ R \in M_{m \times n}(D)$ such that $ {\rm RCS}(R) = W$. }
\end{thm}

\bigskip

\noindent {\bf First proof.} We prove part (i). Part (ii) can be proved similarly. Existence is easy. As $ \dim W \leq m$, we may choose $m$ vectors $ \alpha_1, \ldots, \alpha_m \in D_n$, some of which might be zero,  that span $W$. Set $A \in M_{m \times n}$ to be the matrix whose $i$th row is $ \alpha_i$ for each $ 1 \leq i \leq m$, and hence $ {\rm LRS}(A) = W$. One can easily see that, say by induction on $m$, there exists a left row reduced echelon matrix $R \in M_{m \times n}$   such that $ R =_{lr} A$, which means $ {\rm LRS}(R) = {\rm LRS}(A) = W$, as desired. Now we prove the uniqueness by induction on $n$. Let $ R, R' \in  M_{m \times n}$  be left row reduced echelon matrices with  $ {\rm LRS}(R) ={\rm LRS}(R')= W$. We need to show that $ R = R'$. If  $n=1$, the assertion is easily verified. Suppose the assertion holds for $n-1$. We prove the assertion for $n$. To this end, discard the $n$th column of  $ R $ and $R'$ to obtain row reduced echelon matrices   $R_1 , R'_1 \in M_{m \times (n-1)}$. Clearly, $ {\rm LRS}(R_1) ={\rm LRS}(R'_1)$. So by the inductive hypothesis $ R_1 = R'_1$. Suppose by contradiction that $R \not= R'$. As  $ R_1 = R'_1$, we must have $ {\rm col}_n (R) \not= {\rm col}_n (R')$. Now, let $ X= (x_1 , \ldots, x_n)  \in D^n$ be such that $RX = 0$. This implies that $ R' X = 0$ as well, and hence  $ (R - R') X = 0$. Since  $ R_1 = R'_1$ but $ {\rm col}_n (R) \not= {\rm col}_n (R')$, we see that $ x_n = 0$. Consequently, there are rows of $R$ and $R'$ whose  nonzero leading entries, which is one, occur in the $n$th column of $R$ and $R'$. Clearly, these rows must be the last nonzero rows of $R$ and $R'$ because these leading one entries occur in the $n$th column of $R$ and $R'$. Thus these rows occur in the $r$th row of $R$ and $R'$, where $ r = \dim W$. From this, as $R$ and $R'$ are left row reduced echelon matrices, we see that $ {\rm col}_n (R) = {\rm col}_n (R') = e_r$, where $e_r$ is the column vector with $1$ in the $r$th place and zero elsewhere. This contradicts the hypothesis that   $ {\rm col}_n (R) \not= {\rm col}_n (R')$. Therefore $ R = R'$, which is what we want.

\bigskip 

\noindent {\bf Second proof.} We present a second proof for uniqueness.  Let $ R, R' \in  M_{m \times n}$  be left row reduced echelon matrices with  $ {\rm LRS}(R) ={\rm LRS}(R')= W$. We need to show that $ R = R'$. Let $ k_i$ and $ k'_i$ ($1 \leq i \leq r:= \dim W$) be the column indices of the leading entires of row $i$ of $R$ and $R'$, respectively.  Recall that $(k_i)_{i=1}^r$ and $ (k'_i)_{i=1}^r$ are strictly increasing sequences in $\{ 1, \ldots, n\}$ and that  $ X= (x_1, \ldots, x_n) \in {\rm LRS}(R) ={\rm LRS}(R')= W$ if and only if 
$$X= \sum_{i=1}^r x_{k_i} {\rm row}_i(R) =   \sum_{i=1}^r x_{k'_i} {\rm row}_i(R').$$ 
Thus, it suffices to show that $ k_i = k'_i$ for each $1 \leq i \leq r$. 
We prove this by induction on $i\leq r$. Note that $k_1 = k'_1$ simply because if, for instance,  $ k_1 <k'_1$, then ${\rm row}_1 (R) $ cannot be a linear combination of the rows of $R'$, which is impossible. Thus $ k_1 = k'_1$.  So the assertion holds for $i=1$. Suppose $k_i = k'_i$  for each $ i < i_0 \leq r$. We need to show that $ k_{i_0} = k'_{i_0}$. Again assume, for instance, $ k_{i_0} < k'_{i_0}$.   Then again, ${\rm row}_{i_0} (R) $ cannot be a linear combination of the rows of $R'$, which is impossible.  Therefore,  $ k_i = k'_i$ for each $1 \leq i \leq r$. This completes the proof. \hfill\qed

\bigskip

 In view of the preceding theorem, one can define the operations join and meet on $ \mathcal{R}_{m \times n}$
as follows. Let $R_1 , R_2 \in  \mathcal{R}_{m \times n}$. By $R_1 \vee R_2 $  and $R_1 \wedge R_2 $, we mean the unique matrices in  $\mathcal{R}_{m \times n}$ with the property that $  {\rm LRS} (R_1 \vee R_2)= {\rm LRS} (R_1) + {\rm LRS} (R_2)$  and $  {\rm LRS} (R_1 \wedge R_2)= {\rm LRS} (R_1) \cap {\rm LRS} (R_2)$. For $ C_1 , C_2 \in   \mathcal{C}_{m \times n}$,  one can define $C_1 \vee C_2 $ and $C_1 \wedge C_2 $ in a similar fashion. It is quite straightforward to check that $ (\mathcal{R}_{m \times n},   \vee, \wedge) $ and $ (\mathcal{C}_{m \times n},  \vee, \wedge) $ are modular lattices.   Recall that a  lattice is a triple $(L, \vee, \wedge)$, where $L$ is a nonempty set and $\vee$ and $ \wedge$ are two algebraic operations on $L$ that are commutative, associative, and that they satisfy the absorption laws. Every lattice is a partially ordered set via $  \leq  $, which is naturally defined as follows: $ a \leq b$ if $ a \wedge b = a$, or equivalently $ a \vee b = b$. Consequently,  any isomorphism of  lattices preserves the order structures of them as well. A modular lattice is a lattice that satisfies the modular law, namely, $ a \wedge (b \vee c) = (a \wedge b) \vee (a \wedge c)$ provided that $ b \leq a$ or $ c \leq a$.  It follows from the theorem that  $( \mathcal{R}_{m \times n}, \leq_{lr})$ and  $( \mathcal{C}_{m \times n}, \leq_{rc})$ are partially ordered sets. Note that the partial orders that are induced by the lattice structures of $\mathcal{R}_{m \times n}$ and $\mathcal{C}_{m \times n}$ coincide with $\leq_{lr}$ and $ \leq_{rc}$, respectively.

The set of all left  and right submodules of $ M_{m \times n}(D)$ are, respectively, denoted by $ \mathcal{LS}_{m \times n}(D)$ and $ \mathcal{RS}_{m \times n}(D)$, or simply by $ \mathcal{LS}_{m \times n}$ and   $ \mathcal{RS}_{m \times n}$. By definition, $\mathcal{LS}_{ n}(D) :=   \mathcal{LS}_{n \times n}(D)$ and  $  \mathcal{RS}_{n}(D):= \mathcal{RS}_{n \times n}(D)$. Note that $\mathcal{LS}_{ n}$ and $  \mathcal{RS}_{n}$ are in fact the sets of all left and right ideals of $M_n$, respectively. It is well-known that $\big( \mathcal{LS}_{m \times n}(D), +, \cap \big) $ and  $\big( \mathcal{RS}_{m \times n}(D),  +, \cap\big) $ are modular lattices. It is also well-known that $\big( \mathcal{LS}_{ n}(D), +, . \big) $ and  $\big( \mathcal{RS}_{ n}(D),  +, .\big) $, where $.$ denotes the multiplication of one-sided ideals, are hemirings with left and right identity elements, namely $M_n$, respectively. Recall that a hemiring is a triple $(R, + , .)$, where $R$ is a nonempty set, $(R, +)$ is a commutative monoid with identity element $0$, $(R, .)$ is a semigroup, multiplication distributes over addition from both left and right, and finally $r0 = 0r = 0$ for all $ r \in R$. An element $1_l$ (resp. $1_r$) in a hemiring $R$ is said to be a left (resp. right) identity element  if $ 1_l r = r $ (resp. $r 1_r = r$) for all $ r \in R$.

Let $ r, s \in \mathbb{N}$, $m_i , n_j \in \mathbb{N}$,  $ \sum_{i=1}^r m_i = m$,   $ \sum_{j=1}^s n_j = n$, $M= (m_1 , \ldots, m_r)$, and $N= (n_1,  \ldots, n_s)$. The nest module determined by $ M$ and $N$, denoted by $\mathcal{T}_{(M, N)}(D)$ or simply $\mathcal{T}_{(M, N)}$, is defined as follows
$$  \mathcal{T}_{(M, N)}(D) := \Big\{ (A_{ij}) : 1 \leq i  \leq r,  1 \leq j \leq s, A_{ij} \in M_{m_i \times n_j} (D), A_{ij} = 0 \ \forall \ i > j \Big\}.$$
Naturally, we define  $\mathcal{T}_M(D):=  \mathcal{T}_{(M, M)} (D)$. Again, for the sake of simplicity, we use $ \mathcal{T}_M$ to mean $\mathcal{T}_M(D)$.  It is readily checked that $  \mathcal{T}_{(M, N)}$ is a $(\mathcal{T}_M, \mathcal{T}_N)$-bimodule via the matrix multiplication.
If $ M= N$, then  $\mathcal{T}_M$ is in fact a $(\mathcal{T}_M, \mathcal{T}_M)$-bimodule, or simply a $\mathcal{T}_M$-bimodule, and in particular an $F$-algebra, which we call the nest algebra determined by $M$.  We use $\mathcal{LS}(\mathcal{T}_{(M, N)}(D)) $, or simply $ \mathcal{LS}(\mathcal{T}_{(M, N)}) $, and  $ \mathcal{RS}(\mathcal{T}_{(M, N)}(D))$, or simply $ \mathcal{RS}(\mathcal{T}_{(M, N)})$, to denote  the sets of all left  and right submodules of $  \mathcal{T}_{(M, N)}(D) $, respectively. Again, note that $\mathcal{LS}(\mathcal{T}_{N})$ and $  \mathcal{RS}(\mathcal{T}_{N})$ are in fact the sets of all left and right ideals of 
$\mathcal{T}_N$, respectively.
If $ M= (1, \ldots, 1) \in \mathbb{N}^m$ and $N= (1, \ldots, 1) \in \mathbb{N}^n$,  we set $ \mathcal{T}_{m \times n} :=  \mathcal{T}_{(M, N)} (D)$, which is the set of all upper triangular $m \times n$ matrices. Naturally, we define $ \mathcal{T}_n :=  \mathcal{T}_{n \times n}$, which is the set of all upper triangular matrices of size $n$.

If $A = (A_{ij})$, with   $A_{ij}  \in M_{m_i \times n_j} $ for each  $ 1 \leq i \leq r,  1 \leq j \leq s$, is a block matrix, then we use
$ {\rm Row}_i (A)$ and $ {\rm Col}_j(A) $  to, respectively, denote the $i$th block row and the $j$th block column of $A$. If $A$ is a block matrix with $r$ block rows and $s$ block columns, we use $ A_{ij}$ to denote the $(i, j)$ block entry of $A$.  For a collection $\mathcal F$ of block matrices, one can naturally define $ {\rm Row}_i (\mathcal{F})$, $ {\rm Col}_j(\mathcal{F}) $, and $ \mathcal{F}_{ij}$.  Also for $X  \in M_{m_i \times n_j} $ we use $ \widehat{X}_{ij}$ to denote the block matrix with $X$ in its $(i, j)$ place and zero elsewhere. Now for $ \mathcal{A} \subseteq M_{m_i \times n_j} $, one can define $ \widehat{\mathcal{A}}_{ij}$ in a natural way. A useful observation is in order. Let $ I \in \mathcal{LS}(\mathcal{T}_{(M, N)})$. Then
$$ I =   \left(\begin{array}{ccc}
{\rm Row}_1 (I) \\
\vdots  \\
{\rm Row}_r (I)
\end{array} \right) :=  \left\{  \left(\begin{array}{ccc}
X_1 \\
\vdots  \\
X_r
\end{array} \right):  X_i \in {\rm Row}_i (I)  , \ 1 \leq i \leq r  \right\}.$$
To see this, just note that
$  \left(\begin{array}{ccccc}
0_{m_1 \times n} \\
\vdots\\
0_{m_{i-1}\times n}\\
 {\rm Row}_i(I) \\
 0_{m_{i+1} \times n}\\
\vdots \\
0_{m_r \times n}
\end{array} \right)= {\bf E}_{ii} I \subseteq I $ for each $ 1 \leq i \leq r$, where ${\bf E}_{ii} \in \mathcal{T}_M$ denotes the block matrix with $ I_{m_i}$, the identity matrix of size $m_i$,  in the $(i, i)$ place and zero elsewhere. This clearly implies $ \left(\begin{array}{ccc}
{\rm Row}_1 (I) \\
\vdots  \\
{\rm Row}_r (I)
\end{array} \right) \subseteq I$, for $I$ is additive. The reverse inclusion is trivial. This proves the desired identity.
Likewise, if $J \in  \mathcal{RS}(\mathcal{T}_{(M, N)})$, then
$$ J  =  \left(\begin{array}{ccc}
{\rm Col}_1 (J)  & \cdots &  {\rm Col}_s (J)
\end{array} \right) =  \Big\{  \left(\begin{array}{ccc}
Y_1  & \cdots &  Y_n
\end{array} \right):  Y_j \in {\rm Col}_j (J)  , \ 1 \leq j \leq s  \Big\}.$$
A similar argument establishes the counterparts of the above identities for all $ I \in  \mathcal{LS}_{m \times n}$ and $ J \in \mathcal{RS}_{m \times n}$ in which  $r$, $s$,  ${\rm Row}_i$, and  ${\rm Col}_j$ should be replaced with $m$, $n$,  ${\rm row}_i$, and  ${\rm col}_j$, respectively.

\bigskip

Our first result characterizes one-sided submodules of  $M_{m \times n}$ in terms of left row reduced echelon or right column reduced echelon matrices with entries from $D$.

\bigskip

\begin{thm}\label{1.2}
{\rm  (i)}  {\it There exists an  isomorphism of  lattices
$$ \phi :  \mathcal{LS}_{m \times n} \longrightarrow  \mathcal{R}_{n}$$
 with the property that $ I = M_{m \times n}  \phi(I)$ for all $ I \in  \mathcal{LS}_{m \times n}$. In particular, $\mathcal{LS}_{m \times n}$'s are all isomorphic to $ \mathcal{R}_{n}$ as lattices for all $m \in \mathbb{N}$. }

{\rm (ii)}   {\it There exists an  isomorphism of   lattices
$$ \phi :  \mathcal{RS}_{m \times n} \longrightarrow  \mathcal{C}_{m }$$
 with the property that $ I = \phi(I) M_{m \times n} $ for all $ I \in  \mathcal{RS}_{m \times n}$. In particular, $\mathcal{RS}_{m \times n}$'s are all isomorphic to $ \mathcal{C}_{m}$ as lattices for all $n \in \mathbb{N}$. }
\end{thm}

\bigskip

\noindent {\bf Proof.} We prove part (i). Part (ii) can be proved similarly.  To this end, let $ I \in  \mathcal{LS}_{m \times n}$ be given. For each $ 1 \leq i \leq m$, $ {\rm row}_i (I) \leq D_n$, and hence by Theorem \ref{1.1}, there exists a unique $ R_i \in  \mathcal{R}_{n}$ such that
$ {\rm row}_i (I)= {\rm LRS}(R_i)= D_n R_i$. First, by showing that $ {\rm row}_i (I) =  {\rm row}_1 (I)$, we see that $ R_i = R_1:= R$ for each $ 1 \leq i \leq m$.  As $I$ is a left submodule and the mapping $ A \mapsto E_{i1} A$ takes the $1$st row of $A$ to the $i$th row of it, we have
 $ {\rm row}_1 (I) \subseteq  {\rm row}_i (E_{i1}I) \subseteq  {\rm row}_i (I)$, which obtains  $ {\rm row}_1 (I) \subseteq  {\rm row}_i (I)$. Changing the role of $i$ and $1$,  we obtain $ {\rm row}_i (I) \subseteq  {\rm row}_1 (I)$, and hence $ {\rm row}_i (I) =  {\rm row}_1 (I)$ for each $ 1 \leq i \leq m$. Thus $ {\rm row}_i (I)= D_n R$ for each $ 1 \leq i \leq m$. So we see from the observation we made preceding the theorem that
$$ I =   \left(\begin{array}{ccc}
D_n R \\
\vdots  \\
D_n R
\end{array} \right) = \left\{  \left(\begin{array}{ccc}
X_1 R \\
\vdots  \\
X_m R
\end{array} \right):  X_i \in D_n  , \ 1 \leq i \leq m  \right\}= M_{m \times n}(D) R.$$
Now for $ I  \in  \mathcal{LS}_{m \times n}$, define $\phi(I) = R$. We just showed that $ I = M_{m \times n} R$. Now, let $ I_1 , I_2 \in \mathcal{LS}_{m \times n}$ with $  \phi(I_1) = R_1$ and $  \phi(I_2) = R_2$ be given. Clearly, $ I_1 + I_2 , I_1 \cap I_2 \in \mathcal{LS}_{m \times n}$. On the other hand,
$$ {\rm row}_1 (I_1 + I_2) = {\rm row}_1 (I_1)  +{\rm row}_1 ( I_2) ,   {\rm row}_1(I_1 \cap I_2) =   {\rm row}_1(I_1)  \cap  {\rm row}_1( I_2).$$
 The left equality is easy. We prove the right equality. It is plain that  ${\rm row}_1(I_1 \cap I_2) \subseteq  {\rm row}_1(I_1)  \cap  {\rm row}_1( I_2)$. For the reverse inclusion, let $ X \in {\rm row}_1(I_1)  \cap  {\rm row}_1( I_2)$ be arbitrary. It follows that $ X = {\rm row}_1(A_1) =  {\rm row}_1( A_2)$ for some $ A_1 \in I_1$ and $ A_2 \in I_2$. Then again by the useful observation we made preceding the theorem, we have
$ A:= \left(\begin{array}{cccc}
X \\
0_{ 1\times n}\\
\vdots\\
0_{1 \times n}
\end{array} \right) \in  I_1 \cap I_2 $. Consequently, $ X = {\rm row}_1(A) \in {\rm row}_1(I_1 \cap I_2)$, proving the reverse inclusion, and hence the right equality. Now, from the above equalities, we conclude that  $ \phi(I_1 + I_2) = R_1 \vee R_2$ and $ \phi( I_1 \cap I_2) = R_1 \wedge R_2$. That is, $ \phi$ is a homomorphism of lattices. It remains to show that $ \phi$ is one-to-one and onto. The homomorphism $ \phi$ is one-to-one simply because $ I = M_{m \times n} (D) \phi(I)$ for all $ I \in  \mathcal{LS}_{m \times n}(D)$.  To see that $ \phi$ is onto, let  $ R \in \mathcal{R}_{n}$ be arbitrary. It is quite straightforward to see that $ \phi(I)= R$, where $ I =   M_{m \times n}  R$. This completes the proof.
\hfill\qed

\bigskip

\noindent {\bf Remarks.}  (i) If $ I \in \mathcal{LS}_{m \times n}$ (resp. $ I \in  \mathcal{RS}_{m \times n}$), there exists a unique
$ R \in  \mathcal{R}_{n}$ (resp. $ C \in  \mathcal{C}_{m }$) such that
$$ I = M_{m \times n}  R = \left\{  \left(\begin{array}{ccc}
X_1 R  \\
\vdots  \\
X_m R
\end{array} \right):  X_i \in D_n , \ 1 \leq i \leq m  \right\}$$
(resp.
$$ I =C M_{m \times n}   = \Big\{  \left(\begin{array}{ccc}
CY_1  & \cdots &  CY_n
\end{array} \right):  Y_j \in D^m , \ 1 \leq j \leq n  \Big\}).$$

(ii) The matrix $ R \in  \mathcal{R}_{n}$ (resp. $ C \in  \mathcal{C}_{m }$)  may not be a generator of $I$, for $R$ (resp. $C$) is a square matrix whereas $I$ may consists of rectangular matrices. However, if $m \geq n$ (resp. $m \leq n$), then $I$ is principal and in fact
$ I = M_mR'$ (resp. $ I = C' M_n$), where $R' =   \left(\begin{array}{cc}
R  \\
0_{(m-n) \times n}
\end{array} \right)$ (resp. $C' =  \left(\begin{array}{cc}
C  & 0_{m \times (n -m)}
\end{array} \right)$. Note that if $m=n$, then $ R' = R$ (resp. $C' = C$).

(iii)  If $ I \in  \mathcal{LS}_{m \times n}$ (resp. $ I \in \mathcal{RS}_{m \times n}$), then $ \dim I = m {\rm rank} (R)$ (resp. $ \dim I =  {\rm rank} (C) n$, where $I$ is viewed as a left (resp. right) vector space over $D$.

(iv) In view of the theorem, if $ m=n$, the mapping $\phi$ gives rise to a new operation on $\mathcal{R}_n$ (resp. $\mathcal{C}_n$) as follows. Let $ R_1 , R_2 \in \mathcal{R}_n$ (resp. $C_1 , C_2 \in \mathcal{C}_n$). By definition, $ R_1 \ast R_2 = \phi( M_n R_1. M_n R_2)$ (resp. $ C_1 \ast C_2 = \phi(C_1 M_n . C_2 M_n )$), where $.$ stands for the product of left (resp. right) ideals. It is now clear that $(\mathcal{R}_n, \vee, \ast) $ (resp. $(\mathcal{C}_n, \vee, \ast$))  forms a hemiring with a left (resp. right) identity element and that the mapping $\phi$ is an isomorphism of hemirings for each $n \in \mathbb{N}$. 

\bigskip

Motivated by  \cite[Exercise VIII.3.3]{J1}, we state the following.

\bigskip 

 \begin{cor}\label{1.3}
 {\rm (i)} {\it Let $ I \in  \mathcal{LS}_{ n}$ with $ \phi(I) = R$, where $\phi$ is as in Theorem \ref{1.2}. If $ {\rm rank} (R) = r$, then there exists an invertible matrix $ P \in M_n$ such that $ \phi(P^{-1} I P ) = R'$, where $R' =   \left(\begin{array}{cc}
I_r & 0  \\
0 & 0
\end{array} \right)$.  }

{\rm (ii)} {\it Let $ I \in  \mathcal{RS}_{ n}$ with $ \phi(I) = C$, where $\phi$ is as in Theorem \ref{1.2}. If $ {\rm rank} (C) = r$, then there exists an invertible matrix $ P \in M_n$ such that $ \phi(P^{-1} I P ) = C'$, where $C' =   \left(\begin{array}{cc}
I_r & 0  \\
0 & 0
\end{array} \right)$. }

\end{cor}

\noindent {\bf Proof.} We prove (i). Part (ii) can be proved analogously.  We have $I = M_{ n}  R$. There exists an invertible matrix $P$, which is a product of elementary matrices, such that $R':= RP \in \mathcal{R}_n \cap \mathcal{C}_n$. This clearly implies  $R'=  \left(\begin{array}{cc}
I_r & 0  \\
0 & 0
\end{array} \right)$, where $ r = {\rm rank} (R)$. We can write
$$ P^{-1} I P = P^{-1} M_{ n}  RP = M_{n}R',$$
which implies $ \phi( P^{-1} I P )= R'$, as desired.
 \hfill\qed


  \bigskip

\section{Submodules of nest modules}

In this section we characterize submodules of $\mathcal{T}_{(M, N)}$ in terms of certain finite sequences of left row reduced echelon or right column reduced echelon matrices with entries from $D$. First we need the  following useful lemma.

 \begin{lem}\label{2.1}
 {\rm (i)} {\it  If $ J \in \mathcal{LS}\big(\mathcal{T}_{(M, N)}\big) $, then ${\rm Row}_i (J)  \in \mathcal{LS}_{m_i \times n}$  for all $ 1 \leq i \leq r$.    }

{\rm (ii)} {\it If $ J \in \mathcal{RS}\big(\mathcal{T}_{(M, N)}\big) $, then ${\rm Col}_j (J)  \in \mathcal{RS}_{m \times n_j}$ for all $ 1 \leq j \leq s$.}

\end{lem}

  \bigskip

  \noindent {\bf Proof.}  It suffices to prove (i). Part (ii) can be proved analogously.  Let $ J \in \mathcal{LS}\big(\mathcal{T}_{(M, N)}\big) $ and $ 1 \leq i \leq r$ be given. Clearly, ${\rm Row}_i (J)$ is additive. So it remains to show that $ B {\rm Row}_i (J) \subseteq {\rm Row}_i (J)$ for all $ B \in M_{m_i} $. To see this, given $ B \in M_{m_i}$, we can write
  $$ B {\rm Row}_i (J) =  {\rm Row}_i ({\bf B}_{ii}J)  \subseteq  {\rm Row}_i (J),$$
where   ${\bf B}_{ii} \in \mathcal{T}_M$ is the block matrix with $B$ in the $(i, i)$ place and zero elsewhere.  This completes the proof.
 \hfill\qed

  \bigskip

Let $ N = (n_1, \ldots, n_s)$, where $ n_i \in \mathbb{N}$ and $ n_1 + \cdots + n_s= n$. Define $\mathcal{R}_N(r; D)$ or simply $\mathcal{R}_N(r)$ as follows
$$ \mathcal{R}_N(r ;D)= \mathcal{R}_N(r)$$
$$:= \Big\{ (R_1, \ldots, R_r) \in \mathcal{R}_n^r:   R_1 \geq \cdots \geq R_r, {\rm Col}_j (R_i) = 0\  \forall j < i, 1 \leq i \leq r, 1 \leq j \leq s \Big\}.$$
Define $ \vee$ and $\wedge$ on $\mathcal{R}_N(r)$ componentwise, e.g.,  for $ R = ( R_1, \ldots, R_r)$ and $R'= (R'_1,  \ldots, R'_r)$ in $ \mathcal{R}_N(r; D)$, we write $ R \vee R' := ( R_1 \vee R'_1, \ldots, R_r \vee R'_r)$. Also when $ M = N$, define $\ast$ on $ \mathcal{R}_M (r)$ componentwise via the mapping $\phi$. Note that the operations $ \vee$ and $\wedge$, we just defined on $\mathcal{R}_N(r)$, and the operation $\ast$ defined on $ \mathcal{R}_M (r)$, are well-defined in the sense that  $ R \vee R' ,  R \wedge R'\in \mathcal{R}_N(r)$ whenever $ R , R' \in \mathcal{R}_N(r)$, and that $R \ast R'\in \mathcal{R}_M(r)$ whenever $ R , R' \in \mathcal{R}_M(r)$.  Likewise, for $M= (m_1, \ldots, m_r)$ with $ m_i \in \mathbb{N}$ and $ m_1 + \cdots + m_r= m$, we define
$$  \mathcal{C}_M(s; D) =  \mathcal{C}_M(s)$$
$$ := \Big\{ (C_1, \ldots, C_s) \in \mathcal{C}_m^s:C_1 \leq \cdots \leq C_s, {\rm Row}_i (C_j) = 0\  \forall j < i, 1 \leq i \leq r ,  1 \leq j \leq s\Big\}.$$
Again define $ \vee$, $\wedge$ on $\mathcal{C}_M(s)$ and, when $M=N$, define the operation $\ast$  on $\mathcal{C}_M(r)$ componentwise via $\phi$ in a similar fashion. It is easily checked that $(\mathcal{R}_N(r),  \vee, \wedge)$ and
$(\mathcal{C}_M(s),  \vee, \wedge)$ are modular lattices, and that if $M=N$, then $(\mathcal{R}_M(r),  \vee, \ast)$ and
$(\mathcal{C}_M(r),  \vee, \ast)$ are hemirings. Also, it is readily verified that $\big(\mathcal{LS}(\mathcal{T}_{(M, N)}), +, \cap \big) $ and $ \big( \mathcal{RS}(\mathcal{T}_{(M, N)}),  +, \cap \big )$  are modular lattices and that $\big(\mathcal{LS}(\mathcal{T}_{N}), +, . \big) $ and $ \big( \mathcal{RS}(\mathcal{T}_{N}),  +, . \big )$ are hemirings. Our next result shows that these  modular lattices are pairwise isomorphic via a natural map exhibited in the theorem below.

  \bigskip

 \begin{thm}\label{2.2}
{\rm  (i)}  {\it There exists an  isomorphism of  lattices
$$ \Phi : \mathcal{LS}\big(\mathcal{T}_{(M, N)}\big) \longrightarrow  \mathcal{R}_{N}(r)$$
 with the property that
 $$ I = \left(\begin{array}{cccc}
 M_{m_1 \times n}  R_1\\
  M_{m_2 \times n}  R_2\\
  \vdots \\
   M_{m_r \times n}  R_r\\
 \end{array}\right), $$
 where $ \Phi(I) = R = (R_1 , \ldots , R_r)$ and $R_i =  \phi\big({\rm Row}_i(I)\big)$
  for all $ I \in  \mathcal{LS}\big(\mathcal{T}_{(M, N)}\big)$ and $ 1 \leq i \leq r$.}

{\rm (ii)}   {\it There exists an  isomorphism of   lattices
$$ \Phi :  \mathcal{RS}\big(\mathcal{T}_{(M, N)}\big) \longrightarrow  \mathcal{C}_{M }(s)$$
 with the property that
 $$ I = \left(\begin{array}{cccc}
 C_1 M_{m \times n_1}  &
 C_2  M_{m \times n_2}  &   \cdots  & C_s   M_{m \times n_s}
 \end{array}\right), $$
 where $ \Phi(I) = C = (C_1 , \ldots , C_s)$ with $ C_j =  \phi\big({\rm Col}_j(I)\big)$
 for all $ I \in  \mathcal{RS}\big(\mathcal{T}_{(M, N)}\big)$ and $ 1 \leq j \leq s$.}
\end{thm}

  \bigskip

  \noindent {\bf Proof.} We prove (i). Part (ii) can be proved analogously. Let $ I \in  \mathcal{LS}\big(\mathcal{T}_{(M, N)}\big)$ be given. In view of Lemma \ref{2.1} and the useful observation we made prior to Theorem \ref{1.2}, we see that,  for each $ 1 \leq i \leq r$,  there exists a unique  $ R_i \in  M_n$ such that $ {\rm Row}_i (I) = M_{m_i \times n} R_i$, i.e., $  \phi\big({\rm Row}_i(I)\big)= R_i$, and moreover
$$ I = \left(\begin{array}{cccc}
 M_{m_1 \times n}  R_1\\
  M_{m_2 \times n}  R_2\\
  \vdots \\
   M_{m_r \times n}  R_r\\
 \end{array}\right).$$
It remains to show that $(R_1, \ldots, R_r ) \in    \mathcal{R}_{N}(r)$, and hence  $ \Phi(I) = R = (R_1 , \ldots , R_r)$  is well-defined and that  $ \Phi$ is an isomorphism of modular lattices. First, we show that $ R_q \leq R_p$ whenever $ 1 \leq p, q \leq r$ and $ p < q$. Recall that
$ {\rm Row}_q (I) = M_{m_q \times n}  R_q$. In particular, for a given and arbitrary $ 1 \leq i \leq n$,  we have $E_{1i} R_q \in M_{m_q \times n} R_q= {\rm Row}_q (I)$, where  $E_{1i} \in M_{m_q \times n} $ is a standard matrix. Thus, there exists an $A_i \in I$ such that $E_{1i} R_q  = {\rm Row}_q (A_i)$. Consequently, the $i$th row of $R_q$, which is the first row of $E_{1i} R_q$, occurs as the $i_q$th row of $A_i$, where $ i_q = m_1 + \cdots + m_{q-1} + 1$. Let $ i_p= m_1 + \cdots + m_{p-1} + 1$  if $ p > 1$ and $i_p = 1$ if $p=1$. Clearly, $ i_p < i_q$. Now note that  $ E_{i_p i_q}  \in \mathcal{T}_M$ because $ p < q$, and that $ E_{i_p i_q} A_i$ takes the $i_q$th row of $A_i$, which is in fact the $i$th row of $R_q$, to the $i_p$th row of $A_i$. This means
${\rm row}_i (R_q) = {\rm row}_1 \Big( {\rm Row}_p ( E_{i_p i_q}A_i) \Big)  \in  {\rm row}_1 \Big({\rm Row}_p (I)\Big)$, for $ E_{i_p i_q} \in \mathcal{T}_M$.  That is,  ${\rm row}_i (R_q) \in {\rm row}_1 \Big(M_{m_p \times n}  R_p\Big)$,  which easily implies ${\rm row}_i (R_q) \in {\rm LRS} (R_p)$. But  $ 1 \leq i \leq n$ was arbitrary. Therefore  $ {\rm LRS} (R_q) \leq {\rm LRS} (R_p)$, which means $ R_q \leq R_p$.  Next, we need to show that  ${\rm Col}_j (R_i) = 0$   for all $ j < i $ with $ 1 \leq i \leq r$ and $ 1 \leq j  \leq s$. But this is obvious because for each $ 1 \leq i \leq r$, every row of $R_i$ occurs as a row of ${\rm Row}_i(A) $ for some $ A \in I \subseteq \mathcal{T}_{(M, N)}$.  In view of the useful observation made  regarding the block rows of the left submodules of $\mathcal{T}_{(M, N)}$ prior to Theorem \ref{1.2}, a proof almost identical to that presented in Theorem \ref{1.2} shows that
$$ {\rm Row}_i (I_1 + I_2) = {\rm Row}_i (I_1)  +{\rm Row}_i ( I_2) ,   {\rm Row}_i(I_1 \cap I_2) =   {\rm Row}_i(I_1)  \cap  {\rm Row}_i( I_2),$$
for all  $ I_1 , I_2 \in \mathcal{LS}\big(\mathcal{T}_{(M, N)}\big)$ and $ 1 \leq i \leq r$. This clearly implies that the mapping $\Phi$ is a homomorphism of  lattices. That $\Phi$ is one-to-one follows from the fact that the $\phi$'s are all one-to-one by Theorem \ref{1.2} and that $I_1 = I_2$ if and only if  $ {\rm Row}_i (I_1) = {\rm Row}_i (I_2)$ for each $ 1 \leq i \leq r$. To see that $\Phi$ is onto, for a given $ R = (R_1, \ldots, R_r) \in   \mathcal{R}_{N}(r)$, it is readily checked that $ \Phi(I) = R$, where
$$ I = \left(\begin{array}{cccc}
 M_{m_1 \times n}  R_1\\
  M_{m_2 \times n}  R_2\\
  \vdots \\
   M_{m_r \times n} R_r\\
 \end{array}\right).$$
 This completes the proof.
 \hfill\qed

\bigskip

\noindent {\bf Remarks.} (i) In the special case when $  m= r$,  $ s =n$, $ m_i = 1$,  $n_j = 1$ for all $ 1 \leq i \leq m$ and $ 1 \leq j \leq n$, the theorem characterizes one-sided submodules of upper triangular rectangular matrices in terms of the elements of $\mathcal{R}_{e_n}(m)$
 or $ \mathcal{C}_{e_m}(n)$, where $ e_n \in \mathbb{N}^n$ and $ e_m \in \mathbb{N}^m$ are the elements whose components are all $1$. In the more special case when $m = r = s = n$ and $ m_i = n_j = 1$ for all $ 1 \leq i, j \leq n$,  the theorem characterizes all one-sided ideals of the ring of all upper triangular square matrices with entries from $D$ in terms of the elements of $\mathcal{R}_{e_n}(n)$
 or $ \mathcal{C}_{e_n}(n)$.
 
 (ii) It is clear that the mapping $\Phi$ is an isomorphism of hemirings whenever $ M= N$ and  $\big(\mathcal{LS}(\mathcal{T}_{N}), +, . \big) $, $ \big( \mathcal{RS}(\mathcal{T}_{N}),  +, . \big )$, $(\mathcal{R}_N(r),  \vee, \ast)$, and
$(\mathcal{C}_N(r),  \vee, \ast)$ are viewed as hemirings.  
 
 \bigskip

 \begin{lem}\label{2.3}
 {\rm (i)} {\it  Let $ I \in  \mathcal{LS}_{m \times n}$ with $ \phi(I) = R$. Then $I = M_m A$ for some $ A \in  M_{m \times n} $ if and only if $ {\rm LRS} (R) =  {\rm LRS} (A)$.  }

{\rm (ii)} {\it Let $ I \in  \mathcal{RS}_{m \times n} $ with $ \phi(I) = C$. Then $I =  AM_n$ for some $ A \in M_{m \times n} $ if and only if $ {\rm RCS} (C) =  {\rm RCS} (A)$.}
\end{lem}

  \bigskip

  \noindent {\bf Proof.} We prove (i). Part (ii) can be proved in a similar fashion. First let  $ I \in  \mathcal{LS}_{m \times n} $ with $ \phi(I) = R \in {\rm R}_n$ and $I = M_m A$ for some $ A \in  M_{m \times n} $. We need to show that $ {\rm LRS} (R) =  {\rm LRS} (A)$. By Theorem \ref{1.2},  $I = M_{m\times n} R$. Thus $ I = M_m A=  M_{m\times n} R$. As $ A = I_m A \in I$, we have $ A = X R $ for some $ X \in  M_{m\times n}$. So by the exercise we pointed out preceding Theorem \ref{1.1}, we have  $ {\rm LRS} (A) \subseteq   {\rm LRS} (R)$. To see the reverse inclusion, let $ 1 \leq i \leq n$ be given. It follows that $ E_{1i} R \in I =  M_m(D) A$, where $ E_{1i} \in M_{m\times n}$ is a standard matrix. Thus $ {\rm LRS} (E_{1i} R) \subseteq   {\rm LRS} (A)$. Consequently, $ {\rm row}_i (R) \in {\rm LRS} (A)$. This yields $ {\rm LRS} ( R) \subseteq   {\rm LRS} (A)$, for $ 1 \leq i \leq n$ is arbitrary. Therefore, $ {\rm LRS} (A) =  {\rm LRS} (R)$. Next suppose $ {\rm LRS} (A) =  {\rm LRS} (R)$.  We need to show that  $I = M_mA$. Note that $ A \in I$ because $I = M_{m\times n} R$ and $ A = X R$ for some $ X \in M_{m\times n}$, for $ {\rm LRS} (A) \subseteq  {\rm LRS} (R)$. This yields $  M_m A \subseteq I$. For the reverse inclusion, again since   $ {\rm LRS} (R) \subseteq  {\rm LRS} (A)$, we obtain $ R = P A$ for some $ P \in M_{n \times m}$. But
  $$ I = M_{m\times n} R = M_{m\times n} PA \subseteq M_m A.$$
This completes the proof.
 \hfill\qed

 \bigskip

We need the following proposition to characterize principal submodules of $ \mathcal{T}_{(M,  N)}$.

 \bigskip

 \begin{prop}\label{2.4}
 {\rm (i)} {\it Let $ 0\not= J \in \mathcal{LS}\big(\mathcal{T}_{(M, N)}\big) $ with $ \Phi(J) = (R_1,  \ldots, R_k, 0_n , \ldots, 0_n)$, where $1\leq  k \leq r$ is the largest index for which $ R_k \not= 0$, and $ A \in \mathcal{T}_{(M, N)}$. Then $J =  \mathcal{T}_{M} A$ if and only if
  $${\rm Row}_i (J) =  M_{m_i \times m_i}  {\rm Row}_i (A) + \cdots +  M_{m_i \times m_k }  {\rm Row}_k (A)$$ for all $ 1 \leq i \leq k$ if and only if
  $${\rm LRS}(R_i) =  {\rm LRS}\big( {\rm Row}_i (A)\big) + \cdots + {\rm LRS}\big( {\rm Row}_k (A)\big),$$
   for all $ 1 \leq i \leq k$.
  }

{\rm (ii)} {\it Let $ 0\not= J \in \mathcal{RS}\big(\mathcal{T}_{(M, N)}\big) $  with $ \Phi(J) = (  0_m, \ldots, 0_m, C_k, \ldots, C_s)$, where $ 1 \leq k \leq s$ is the smallest index for which $ C_k \not= 0$, and $ A \in \mathcal{T}_{(M, N)}$. Then $J = A \mathcal{T}_{N} $ if and only if
 $$ {\rm Col}_j(J) =  {\rm Col}_k (A)  M_{n_k \times n_j }   + \cdots +  {\rm Col}_j (A)M_{ n_j \times n_j } $$
  for all $ k \leq j \leq s$ if and only if
  $${\rm RCS}(C_j) = {\rm RCS}\big( {\rm Col}_k (A) \big) + \cdots +  {\rm RCS}\big( {\rm Col}_j (A) \big)$$
  for all $ k \leq j \leq s$.
 }
\end{prop}

  \bigskip

  \noindent {\bf Proof.} The assertion easily follows from the useful observations we have already made, namely  $J =  \mathcal{T}_{M} A$  if and only if  $M_{m_i \times n } R_i = {\rm Row}_i (J) = {\rm Row}_i (\mathcal{T}_{M} A)$ for all $ 1 \leq i \leq r$, and that $J = A \mathcal{T}_{N} $ if and only if   $C_j M_{m \times n_j }=  {\rm Col}_j(J) =  {\rm Col}_j (A \mathcal{T}_{N})$
  for all $ 1 \leq j \leq s$.
 \hfill\qed

  \bigskip

 \begin{cor}\label{2.5}
  {\it Let  $ A \in  \mathcal{T}_{m \times n}$.  Then  }

\noindent  {\rm (i)} {\it given  $ J \in  \mathcal{LS}( \mathcal{T}_{m \times n} )$,  $ J =\mathcal{T}_{m} A$ if and only if   $ {\rm row}_i (J)= \langle \{{\rm row}_k (A): i \leq k \leq m \}\rangle $ for each $ 1 \leq i \leq m$. }

\noindent {\rm (ii)} {\it given  $ J \in  \mathcal{RS}( \mathcal{T}_{m \times n} )$,  $ J =A \mathcal{T}_{n} $ if and only if $ {\rm col}_j (J)= \langle \{{\rm col}_k (A): 1 \leq k \leq j \}\rangle $ for each $ 1 \leq j \leq n$. }
\end{cor}
 \bigskip

  \noindent {\bf Proof.} This corollary is a special case of the preceding proposition.
 \hfill\qed

 \bigskip

The following characterizes the principal submodules of  $ \mathcal{T}_{(M,  N)}$.

\bigskip

 \begin{thm}\label{2.6}
 {\rm (i)} {\it Let  $0\not=  J \in  \mathcal{LS}( \mathcal{T}_{(M, N)} )$ with $ \Phi(J) = ( R_1, \ldots, R_k, 0_n,  \ldots, 0_n)$, where $1\leq  k \leq r$ is the largest index for which $ R_k \not= 0_n$. Then $J$ is principal if and only if  $ {\rm rank}(R_{i}) - {\rm rank}(R_{i+1}) \leq m_i$ for each $ 1 \leq i \leq k$. By convention, $ R_{r+1}:= 0_n$.}

{\rm (ii)} {\it  Let  $ 0\not= J \in  \mathcal{RS}( \mathcal{T}_{(M, N)} )$ with $ \Phi(J) = (  0_m,  \ldots, 0_m, C_k, \ldots, C_s)$, where $ 1 \leq k \leq s$ is the smallest index for which $ C_k \not= 0_m$. Then $J$ is principal if and only if  $ {\rm rank}(C_{j}) - {\rm rank}(C_{j-1}) \leq n_j$ for each $ k \leq j \leq s$. By convention, $ C_{0}:= 0_m$.}

\end{thm}

  \bigskip

  \noindent {\bf Proof.} We prove (i). Part (ii) can be proved analogously. First, let  $ J = \mathcal{T}_{M} A$ for some  $ A \in \mathcal{T}_{(M, N)}$. It follows from Proposition \ref{2.4} that
  $${\rm LRS}(R_i) =  {\rm LRS}\big( {\rm Row}_i (A)\big) + \cdots + {\rm LRS}\big( {\rm Row}_k (A)\big).$$
It is plain that
\begin{eqnarray*}
 {\rm rank}(R_{i}) & =&  \dim {\rm LRS}(R_i) \\
  & \leq&  \dim {\rm LRS}\big( {\rm Row}_i (A)\big) + \dim\Big(  {\rm LRS}\big( {\rm Row}_{i+1} (A)\big) + \cdots + {\rm LRS}\big( {\rm Row}_k (A)\big) \Big)\\
  & \leq  & m_i +  \dim {\rm LRS}(R_{i+1}),
\end{eqnarray*}
and hence   $ {\rm rank}(R_{i}) - {\rm rank}(R_{i+1}) \leq m_i$ for each $ 1 \leq i \leq k$, as desired. Next, let $ {\rm rank}(R_{i}) - {\rm rank}(R_{i+1}) \leq m_i$ for each $ 1 \leq i \leq k$. In view of Proposition \ref{2.4}, it suffices to find an $A \in \mathcal{T}_{(M, N)}$ such that
$${\rm LRS}(R_i) =  {\rm LRS}\big( {\rm Row}_i (A)\big) + \cdots + {\rm LRS}\big( {\rm Row}_k (A)\big),$$
   for all $ 1 \leq i \leq k$. We find $A$ by finding  its block rows, i.e., ${\rm Row}_i (A)$'s where $ 1 \leq i \leq r$. To this end, set  ${\rm Row}_i (A) = 0_{m_i \times n}$ for each $ k < i \leq r$. By hypothesis  ${\rm rank} (R_k) \leq m_k$. Thus there exists an $A_k \in  M_{m_k \times n}$, whose rows are chosen from those of  $ {\rm LRS}(R_k)$ or are zero, such that $ {\rm LRS}(A_k) ={\rm LRS}(R_k)$. Set $  {\rm Row}_k (A) =A_k$. Now as $ {\rm rank} (R_{k-1})  - {\rm rank} (R_k) \leq m_{k-1}$ and $ {\rm LRS}(R_k) \subseteq  {\rm LRS}(R_{k-1}) $,
we can find    $A_{k-1} \in  M_{m_{k-1} \times n}$,  whose rows are chosen from those of  ${\rm LRS}( R_{k-1})$ or are zero,  such that $ {\rm rank}(A_{k-1}) =  {\rm rank} (R_{k-1})  - {\rm rank} (R_k) \leq m_{k-1}$. Continuing in this way, we obtain an $A \in  \mathcal{T}_{(M, N)}$ such that
$$  {\rm rank}({\rm Row}_i (A))  ={\rm rank}(A_i)=  {\rm rank} (R_{i-1})  - {\rm rank} (R_i) \leq m_i $$
for all $ 1 \leq i \leq k$. It remains to show that  $${\rm LRS}(R_i) =  {\rm LRS}\big( {\rm Row}_i (A)\big) + \cdots + {\rm LRS}\big( {\rm Row}_k (A)\big),$$
   for all $ 1 \leq i \leq k$. Fix $ 1 \leq i \leq k$ and let $ \alpha $ be an arbitrary nonzero row of  $R_i$. As $ R_i \geq \cdots \geq R_k$ there exists a largest $i \leq i' \leq k $ such that  $ \alpha $ is in ${\rm LRS}(R_i')$, and hence  not in ${\rm LRS}(R_{i'+1})$.  But this means $ \alpha$ is in ${\rm LRS}(A_{i'})= {\rm LRS}\big({\rm Row}_{i'} (A)\big)$. Consequently,
 $$ \alpha \in {\rm LRS}\big( {\rm Row}_{i'} (A)\big) \subseteq {\rm LRS}\big( {\rm Row}_i (A)\big) + \cdots + {\rm LRS}\big( {\rm Row}_k (A)\big).$$
 This implies
 $${\rm LRS}(R_i) \subseteq  {\rm LRS}\big( {\rm Row}_i (A)\big) + \cdots + {\rm LRS}\big( {\rm Row}_k (A)\big),$$
 for $ \alpha $ was an arbitrary nonzero row of  $R_i$. The reverse inclusion is trivial because by the way we constructed $A_i$'s, we have
  $$  {\rm LRS}\big( {\rm Row}_p (A)\big) \subseteq {\rm LRS}(R_p)  \subseteq {\rm LRS}(R_i)$$
  for each $ i \leq p \leq k$. This yields
  $${\rm LRS}\big( {\rm Row}_i (A)\big) + \cdots + {\rm LRS}\big( {\rm Row}_k (A)\big) \subseteq {\rm LRS}(R_i),$$
  which is what we want. This completes the proof.
 \hfill\qed

\bigskip

  \bigskip

 \begin{cor}\label{2.7}
 {\rm (i)} {\it Let  $0\not=  J \in  \mathcal{LS}( \mathcal{T}_{m \times n} )$ with $ \Phi(J) = ( R_1, \ldots, R_k, 0_n \ldots, 0_n)$, where $1\leq  k \leq m$ is the largest index for which $ R_k \not= 0_n$. Then $J$ is principal if and only if  $ {\rm rank}(R_{i}) - {\rm rank}(R_{i+1}) \in \{0,  1\}$ for each $ 1 \leq i \leq k$. By convention, $ R_{r+1}:= 0_n$.}

{\rm (ii)} {\it  Let  $0\not=  J \in  \mathcal{RS}( \mathcal{T}_{m \times n} )$ with $ \Phi(J) = (  0_m, \ldots, 0_m, C_k, \ldots, C_n)$, where $ 1 \leq k \leq n$ is the smallest index for which $ C_k \not= 0_m$. Then $J$ is principal if and only if  $ {\rm rank}(C_{j}) - {\rm rank}(C_{j-1}) \in \{0,  1\}$ for each $ k \leq j \leq n$. By convention, $ C_{0}:= 0_m$.}

\end{cor}

  \bigskip

  \noindent {\bf Proof.} This is a quick consequence of the preceding theorem.
 \hfill\qed

\bigskip

\section{Subbimodules of nest modules}

\bigskip

We now intend to characterize subbimodules of nest modules. As it turns out the only subbimodules of $M_{m \times n}$ are $ 0$ and  $M_{m \times n}$ itself. A proof can be given based on  \cite[Corollary 1.10]{RY2}. Here, we present a proof based on one of the facts we left as an exercise preceding Theorem \ref{1.1}.

\bigskip

 \begin{lem}\label{3.1}
 {\it The only subbimodules of $M_{m \times n}$ are the trivial ones, namely $ 0$ and $M_{m \times n}$.  }
\end{lem}

  \bigskip

  \noindent {\bf Proof.} Let $I$ be a nonzero subbimodule of $M_{m \times n}$. We show that $ I =  M_{m \times n}$. As $ I \not= 0$, there is a nonzero element $A=(a_{ij}) \in I$.   Thus there are $ 1 \leq i_0 \leq m $ and $ 1 \leq j_0 \leq n$ such that $ a_{i_0 j_0} \not= 0$. Let $ B= (b_{ij})  \in M_{m \times n}$ be arbitrary. Set $ K := \{(i, j ) : 1 \leq i \leq m, 1 \leq j \leq n, b_{ij} \not= 0\}$. We can write
  $ B = \sum_{(i, j) \in K} B_{ij}$, where $ B_{ij} = E_{ii} B E_{jj}$, which is the matrix whose entries are all zero except for its $(i, j)$ entry which is the same as that of $B$. Note that $ E_{ii} \in M_m $ and $ E_{jj} \in M_n $ are standard matrices. Clearly,
  $${\rm LRS}(B_{ij}) \leq \langle e_j  \rangle = {\rm LRS}( E_{ii_0} A  E_{jj_0}),$$
  where $ e_j \in D_n$ is the row vector whose components are all zero except for its $i$th component which is $1$.
  Again $ E_{ii_0} \in M_m $ and $ E_{jj_0} \in M_n $ are standard matrices.
  Thus there is a  $ C_{ij} \in M_m$ such that $ B_{ij} = C_{ij}  E_{ii_0} A  E_{jj_0}$, implying that $  B_{ij} \in I$ for each $(i, j) \in K$. This yields $ B \in I$. Therefore, $ I = M_{m \times n}$, as desired.
 \hfill\qed

\bigskip

\noindent {\bf Remark.} It is possible to use Theorem \ref{1.2}, to present a short proof of this lemma. To this end, viewing the given nonzero subbimodule $I$ as a left submodule, it suffices to show that ${\rm rank}(R) = n$, where $ R = \phi(I)$. As $I$ is a nonzero subbimodule of $M_{m \times n}$, we have $ IP = I$ for any invertible $ P \in M_n$. So we must have $R \in \mathcal{R}_n \cap \mathcal{C}_n$, implying that  $R =   \left(\begin{array}{cc}
I_r & 0  \\
0 & 0
\end{array} \right)$, where $ r = {\rm rank}(R)$. But $r = n$, for otherwise multiplying $I$ from the right by the permutation matrix $P$ obtained by exchanging the $r$th column and the $n$th column of the identity matrix, we see that $ \phi(I) $ could be both $R, RP \in  \mathcal{R}_n$, which is impossible as $ R \not= RP$. 

\bigskip

The following characterizes all subbimodules of nest modules of matrices.

  \bigskip

 \begin{thm}\label{3.2}
 {\it Let $J$ be a nonzero subbimodule of $ \mathcal{T}_{(M, N)}$. Then there exists a unique $ 1 \leq i \leq \min (r, s)$, which depends on $J$,  a unique increasing sequence $(j_1, \ldots, j_i)$ such that  $k \leq j_k \leq s$ and $ j_k < j_{k+1} $ if $ j_k = k< i$ for all $ 1 \leq k \leq i$,  and a unique sequence $(R_1 , \ldots, R_i)$ with
 $ R_1 = I_n$ if $j_1 = 1$, $R_1 =\left(\begin{array}{ll}
0_{r_1 \times (n-r_1)} &  I_{r_1} \\
0_{n-r_1} & 0_{(n-r_1) \times r_1}
\end{array} \right) \in M_n $ if $ j_1 > 1$, and $R_k =\left(\begin{array}{ll}
0_{r_k \times (n-r_k)} &  I_{r_k} \\
0_{n-r_k} & 0_{(n-r_k) \times r_k}
\end{array} \right) \in M_n$ if $ k\geq 2$, where
$ r_k = \sum_{l= j_k}^s n_l$ provided $j_k > 1$, such that
$$ J = \left(\begin{array}{ccccc}
 M_{m_1 \times n}  R_1\\
   \vdots \\
  M_{m_i \times n}  R_i\\
  0_{m_{i+1} \times n} \\
  \vdots \\
 0_{m_r \times n}
 \end{array}\right).$$
Moreover, every subbimodule of $ \mathcal{T}_{(M, N)}$ is principal.  }

\end{thm}

  \bigskip

  \noindent {\bf Proof.} Let $J$ be a nonzero subbimodule of $ \mathcal{T}_{(M, N)}$ and $1 \leq  i \leq r$ be the largest integer for which $ {\rm Row}_i (J) \not= 0$. If $ r \leq s$, then $ i \leq r = \min (r, s)$. If $ r > s$, then $ {\rm Row}_k (J) = 0 $ for each $  k > r$ because $ I \subseteq   \mathcal{T}_{(M, N)}$. Thus $ i \leq s= \min(r, s)$. Now we show that $J_{kl}$ is a subbimodule of $M_{m_k \times n_l}$ for each $ 1 \leq k \leq r$ and $ 1 \leq l \leq s$.   Let $ X \in J_{kl}$ and $ B \in M_{m_k}$ and $ C \in M_{n_l}$ be arbitrary. It follows that $ X = A_{kl}$ for some $ A \in J$.  We can write $ \widehat{B}_{kk} A   \widehat{C}_{ll} =  \widehat{BXC}_{kl}$. But $ A \in J$ and $ J$ is a subbimodule. Thus $BXC \in  J_{kl}$, and hence  $J_{kl}$ is a subbimodule of $M_{m_k \times n_l}$. It follows from the preceding lemma that $  J_{kl} = 0_{m_k \times n_l}$ or $  J_{kl} =M_{m_k \times n_l}$. For $ 1 \leq k \leq i$, let $ k \leq j_k \leq s$ be the smallest integer  for which $ J_{kj_k} \not= 0_{m_k \times n_{j_k}}$, or equivalently, $ J_{kj_k}= M_{m_k \times n_{j_k}}$. So there exists an $A \in J$ such that $A _{kj_k} = E_{11} \in   M_{m_k \times n_{j_k}}$.  For $ 1 \leq u \leq i$, let $ p_u =1 +  \sum_{t=1}^{u-1}m_t $.    Clearly, $ E_{p_k p_l}  \in M_{m \times n} \cap \mathcal{T}_M$ for each $ 1\leq k < l < i$. Thus $  E_{p_k p_l} A \in J$. On the other hand,  the $p_k$th row of $  E_{p_k p_l} A \in J$ is the same as its $ p_l$th row, and hence $A _{kj_l} \not= 0$. Therefore,  $ J_{kj_l} \not= 0_{m_k \times n_{j_l}}$. This yields $ j_k \leq  j_l$. That is, the sequence $(j_1, \ldots, j_i)$ is increasing. Note that if $ j_k = k< i$, then $ j_k < j_{k+1} $, for $ J \subseteq \mathcal{T}_{(M, N)}$. To prove the equality, note first that
   $$ J = \left(\begin{array}{ccccc}
{\rm Row}_1 (J)\\
   \vdots \\
{\rm Row}_i (J)\\
  0_{m_{i+1} \times n} \\
  \vdots \\
 0_{m_r \times n}
 \end{array}\right), $$
for $ {\rm Row}_k(J) = 0_{m_k \times n}$ for each $ k > i$.  But $ {\rm Row}_k (J)= \left(\begin{array}{ccc} 0_{m_k \times x_k} & M_{m_k \times y_k}   \end{array}\right)= M_{m_1 \times n} (D) R_k$, where $ x_k = \sum_{j=1}^{j_{k}-1} n_j$ and $y_k =  \sum_{j=j_k}^{s} n_j$. Note that if $ j_k = 1$, then $ k= 1$, $x_1= 0$, $y_1= n$, and $R_1= I_n$, and hence $ {\rm Row}_1 (J)=  M_{m_1 \times n} = M_{m_1 \times n} (D) R_1$. Conversely,  with an  $ 1 \leq i \leq \min (r, s)$, an increasing sequence $(j_1, \ldots, j_i)$ with $k \leq j_k \leq s$ and $ j_k < j_{k+1} $ if $ j_k = k< i$ for all $ 1 \leq k \leq i$, and  a sequence $(R_1 , \ldots, R_i)$  as in the statement of the theorem, a simple calculation with block matrices shows that
$$ \left(\begin{array}{ccccc}
 M_{m_1 \times n} R_1\\
   \vdots \\
  M_{m_i \times n}  R_i\\
  0_{m_{i+1} \times n} \\
  \vdots \\
 0_{m_r \times n}
 \end{array}\right)$$
 is in fact a subbimodule of    $ \mathcal{T}_{(M, N)}$. Clearly the subbimodule $\{0_{m \times n}\}$ is principal. So let  $$ J = \left(\begin{array}{ccccc}
 M_{m_1 \times n}  R_1\\
   \vdots \\
  M_{m_i \times n}  R_i\\
  0_{m_{i+1} \times n} \\
  \vdots \\
 0_{m_r \times n}
 \end{array}\right)$$
be a nonzero subbimodule of $ \mathcal{T}_{(M, N)}$.   Then $J$ is generated by any matrix $ A \in J$ for which $ A_{kj_k} \not= 0$ for all $ 1 \leq k \leq i$, e.g., for
 $$ A = \left(\begin{array}{ccccc}
A_1\\
   \vdots \\
A_i\\
  0_{m_{i+1} \times n} \\
  \vdots \\
 0_{m_r \times n}
 \end{array}\right), $$
 where  $A_k = \left(\begin{array}{ccc} 0_{m_k \times x_k} & E_{m_k 1}    \end{array}\right) \in M_{m_1 \times n}  R_k$, where $  E_{m_k 1} \in M_{m_k \times y_k}$ for each $ 1 \leq k \leq i$. Clearly, the subbimodule generated by $A$ shares the same $ 1 \leq i \leq \min (r, s)$ and $j_k$'s  ($ 1 \leq k \leq i$)  as $J$. This implies that $J$ is generated by $A$. So the proof is complete.
 \hfill\qed

\bigskip

\noindent {\bf Remarks.}  (i) In the special case $M=N$, the theorem characterizes two-sided ideals of the nest algebra $ \mathcal{T}_{M}$. Moreover any two-sided ideal of $ \mathcal{T}_{M}$ is principal. 

(ii) In the special case when $  m= r$,  $ s =n$, $ m_i = 1$,  $n_j = 1$ for all $ 1 \leq i \leq m$ and $ 1 \leq j \leq n$, the theorem characterizes subbimodules of upper triangular rectangular matrices and in the more special case when $m = r = s = n$ and $ m_i = n_j = 1$ for all $ 1 \leq i, j \leq n$,  the theorem characterizes all two-sided ideals of the ring of all upper triangular square matrices with entries from $D$.

 \bigskip

\noindent  {\bf Acknowledgment.} This paper was submitted while the second named author was on sabbatical leave at the University of Waterloo. He would like to thank Department of Pure Mathematics of the University of Waterloo, and in particular Professors Laurent Marcoux and Heydar Radjavi, for their support.

\end{document}